\documentclass[10pt]{amsart}
\usepackage{amssymb}
\usepackage{amscd}

\def\today{\number\day\space\ifcase\month\or   January\or February\or
   March\or April\or May\or June\or   July\or August\or September\or
   October\or November\or December\fi\   \number\year}

\theoremstyle{definition}
\newtheorem{thm}{Theorem}[section]
\newtheorem{lem}[thm]{Lemma}
\newtheorem{prp}[thm]{Proposition}
\newtheorem{dfn}[thm]{Definition}
\newtheorem{cor}[thm]{Corollary}

\newtheorem{rmk}[thm]{Remark}

\newtheorem{exa}[thm]{Example}

\newcommand{\beq}{\begin{equation}}
\newcommand{\eeq}{\end{equation}}
\newcommand{\beqr}{\begin{eqnarray*}}
\newcommand{\eeqr}{\end{eqnarray*}}
\newcommand{\bal}{\begin{align*}}
\newcommand{\eal}{\end{align*}}
\newcommand{\bei}{\begin{itemize}}
\newcommand{\eei}{\end{itemize}}

\newcommand{\dist}{{\mathrm{dist}}}

\newcommand{\ep}{\varepsilon}

\newcommand{\ph}{\varphi}

\newcommand{\Z}{{\mathbb{Z}}}

\newcommand{\C}{{\mathbb{C}}}
\newcommand{\N}{{\mathbb{N}}}

\newcommand{\ca}{C*-algebra}

\newcommand{\Index}{{\mathrm{Index}}}
\newcommand{\id}{{\mathrm{id}}}

\title[Strongly self-absorbing property for inclusions of $C^*$-albebras] 
{Strongly self-absorbing property for inclusions of $C^*$-algebras with a finite Watatani index}
\author{Hiroyuki Osaka$^*$}
\date{Feb. 22, 2010}
\thanks{$^*$Research of the first author partially supported by the JSPS grant for Scientific Research No.20540220}
\address{ Department of Mathematical Sciences\\
  Ritsumeikan University\\ Kusatsu, Shiga, 525-8577  Japan}

\email[]{osaka@se.ritsumei.ac.jp}

\author{Tamotsu Teruya}

\address{Department of Mathematical Sciences\\
  Ritsumeikan University\\ Kusatsu, Shiga, 525-8577  Japan}

\email[]{teruya@se.ritsumei.ac.jp}


\begin{document}
\maketitle

\begin{abstract}
Let $P \subset A$ be a inclusion of unital C*-algebras and 
$E\colon A \rightarrow P$ be a conditional expectation of index finite type.
We introduce a Rokhlin property  for $E$ and discuss about $\mathcal{D}$-absorbing proeprty,
where $\mathcal{D}$ is a separable, unital, strongly self-absorbing C*-algebra, i.e., 
$\mathcal{D} \not\cong \C$ and there is an isomorphism 
$\varphi\colon \mathcal{D} \rightarrow \mathcal{D} \otimes \mathcal{D}$ 
such that $\varphi$ is approximate unitarily equivalent to the embedding $d \mapsto d \otimes 1_D$.
an UHF algebra of infinite type, Jiang-Su algebra $\mathcal{Z}$, and Cuntz algebras $O_2$, $O_\infty$ are
typical examples for strongly self-absorbing C*-algebras.
In this paper we consider permanent properties for strongly self-absorbing property under inclusions 
of unital C*-algebras with a finite Watatani index. We show the followings:

Let $P \subset A$ be an inclusion of unital C*-algebras, $E$ a conditional 
expectation from $A$ onto with a finite index, and $\mathcal{D}$ be a separable unital 
strongly self-absorbing C*-algebra.
\begin{enumerate}
\item
 If $A$ is a separable $\mathcal{D}$-absorbing 
(i.e. $A \otimes \mathcal{D} \cong A$), and $E$ has the Rokhlin property, then $P$ is 
$\mathcal{D}$-absorbing.
\item
If $\mathcal{D}$ is the universal UHF algebra $\mathcal{U}_\infty$ and 
$\alpha\colon G \rightarrow Aut(D)$ has the Rokhlin property, 
then $\mathcal{U}_\infty \rtimes_\alpha G \cong \mathcal{U}_\infty$.
\item
If $A = \mathcal{D}$ is semiprojective and $E$ has the Rokhlin property, then 
$P \cong \mathcal{D}$. 
In particlular this is true when $\mathcal{D}$ is either $\mathcal{O}_2$ or 
$\mathcal{O}_\infty$.
\end{enumerate}
As an application, under the condition that  
$A$ is a unital $\mathcal{D}$-absorbing C*-algebra
and $\alpha$ is an action from a finite group $G$ on $A$ 
with the Rokhlin property we have then

\begin{enumerate}
\item
the crossed product algebra 
$A \rtimes_\alpha G$ is also $\mathcal{D}$-absorbing. 
\item for any subgroup $H \subset G$ $A^H$ is $\mathcal{D}$-absorbing.
\item
if $A$ is $\mathcal{O}_2$, then 
$A \rtimes_\alpha G \cong \mathcal{O}_2$.
\end{enumerate}
\end{abstract}

\section{Introduction}
A separable C*-algebra $D$ is said to be strongly self-absorbing if $D \not= \C$ and 
there exists an isomorphism $\varphi\colon D \rightarrow D \otimes D$ such that $\varphi$ and 
$id_D \otimes 1_D$ are approximate unitarily equivalent *-homomorphism.
Strongly self-absorbing C*-algebras are automatically simple, nuclear, and have at most 
one tracial state \cite{KP:embedding}. 
The UHF algebra of infinite type and the 
Jiang-Su algebra $\mathcal{Z}$, and the Cuntz algebras $\mathcal{O}_2$ and $\mathcal{O}_\infty$
are strongly self-absorbing.

In Elliott program to classify nuclear C*-algebras by K-theory data
the systematic use of strongly self-absorbing C*-algebras plays a central role. 
In the purely infinite case the Cuntz algebra $\mathcal{O}_\infty$ is a cornerstone of 
the Kirchberg- Phillips classification of simple purely infinite C*-algebras 
\cite{KP:embedding} \cite{Phillips:classification}.
In the stably finite case the Jiang-Su algebra $\mathcal{Z}$ plays a role 
similar  to that of $\mathcal{O}_\infty$.
In fact Jiang-Su proved in \cite{Jiang-Su:absorbing} 
 that simple, infinite dimensional AF algebras and 
Kirchberg algebras (simple, nuclear, purely infinite and satisfying the Universal Coefficient Theorem) are $\mathcal{Z}$-stable, that is, for any such an algebra $A$ 
one has an isomorphism $\alpha\colon A \rightarrow A \otimes \mathcal{Z}$.
Gong, Jiang, and Su proved in \cite{GJS} that $(K_0(A), K_0(A)^+)$ 
is isomorphic to $(K_0(A \otimes \mathcal{Z}), K_0(A \otimes \mathcal{Z})_+)$ 
if and only if $K_0(A)$ is weakly unperforated as an ordered group, when $A$ is a simple 
C*-algebra. Hence $A$ and $A \otimes \mathcal{Z}$ have isomorphic Elliott invariant if 
$A$ is simple with weakly unperforated $K_0$-group, that is, $A \cong A \otimes \mathcal{Z}$
whenever $A$ is classifiable. On the contrary, R\o rdam and Toms 
in \cite{R: finite and infinite} and \cite{T:K-theory} 
presented examples which have the same Elliott invariant as, but are not isomorphic to, and 
not $\mathcal{Z}$-absorbing.  
So it appears plausible that the Elliott conjecture, which is formulated in \cite{Ro}, 
holds for all simple, unital, nuclear, separable $\mathcal{Z}$-absorbing C*-algebras.

In this paper we begin by reconsidering the $\mathcal{D}$-absorbing property 
for crossed product of a C*-algebra $A$ with $\mathcal{D}$-absorbing by a finite group action 
with the Rokhlin property in the framework of inclusion of unital C*-algebras $P \subset A$
of Watatani index finite (\cite{Watatani:index}) and show that 
if a faithfule conditional expectation $E$ from $A$ to $P$ has the Rokhlin property in the sense of \cite{KOT}, then
$P$ is $\mathcal{D}$-absorbing. 
In section 4 we study the  strongly self-absorbing property and 
we show that the universal UHF algebra $\mathcal{U}_\infty$ 
is stable under crossed product of a finite group action with the Rokhlin 
property. 
In the case of pure infiniteness we show that 
if $P \subset \mathcal{O}_2$ (resp. $P \subset O_\infty$) is an inclusion of 
unital C*-algebra with the Rokhlin property, then $P \cong \mathcal{O}_2$ (resp. $P \cong O_\infty$). 
In section 5 we prove that the $\mathcal{D}$-absorbing property is stable 
under local $\mathcal{C}$-algebra which was defined in \cite{OP:Rohlin}. 
Let $\mathcal{C}$ be the class of all separable, unital, 
$\mathcal{D}$- absorbing C*-algebras. We prove that if $A$ is locally approximated by 
some C*-algebras in $\mathcal{C}$, then $A \in \mathcal{C}$.
We  stress that the strongly self-absorbing is not generarily stable under crossed products 
by actions with the tracial Rokhlin property, that is, there is a symmetry with the tracial Rokhlin property 
on $\mathcal{U}_\infty$ such that $\mathcal{U}_\infty \rtimes_\alpha \Z/2\Z$ is not strongly self-absorbing. 
In the last section we present examples of inclusion of unital C*-algebras with 
the Rokhlin property which do not come from finite group actions.

The first auhtor would like to thank Hiroki Matui for a fruitful discussion.

\section{Rokhlin property}
A notion of the Rokhlin property has already appeared in \cite{HJ:82}, \cite{HJ:83}, 
and \cite{Kishimoto:77} with a different name. Indeed, Herman and Jones investigated in 
\cite{HJ:82} and \cite{HJ:83} a class of finite group actions with what we call here 
the Rokhlin property on UHF C*-algebra. 
Izumi introduced  the term "the Rokhlin property" for a finite group action 
in \cite[Definition 3.1]{Izumi:Rohlin1} as follows and characterized two 
finite group action with the Rokhlin property on a separable unital C*-algebra 
using the Evans-Kishimoto intertwining argument in \cite{EK:intertwining}.

\vskip 3mm

For a $C^*$-algebra $A$, we set 
\begin{eqnarray*}
c_0(A) &=& \{(a_n) \in l^\infty(\N, A): \lim\limits_{n \to \infty} \Vert a_n \Vert= 0\}  \\
 A^\infty &=&l^\infty(\N, A)/c_0(A).  
\end{eqnarray*}
We identify $A$ with the $C^*$-subalgebra of $A^\infty$ consisting of the equivalence classes of constant sequences and 
set
$$
A_\infty = A^\infty \cap A'.
$$
For an automorphism $\alpha \in {\rm Aut}(A)$, we denote by $\alpha^\infty$ and $\alpha_\infty$ the automorphisms of 
$A^\infty$ and $A_\infty$ induced by $\alpha$, respectively.

\vskip 3mm

\begin{dfn}\label{def:group action}
Let $\alpha$ be an action of a finite group $G$ on a unital $C^*$-algebra $A$. 
$\alpha$ is said to have the {\it Rokhlin property} if there exists a partition of unity 
$\{e_g\}_{g \in G} \subset A_\infty$ consisting of projections satisfying 
$$
(\alpha_g)_\infty(e_h) = e_{gh} \quad \text{for} \  g, h \in G.
$$
We call $\{e_g\}_{g\in G}$  Rokhlin projections. 
\end{dfn}  

\vskip 3mm

Motivated by Definition~\ref{def:group action}
Kodaka, Osaka, and Teruya introduced the Rokhlin property for a inclusion of 
unital C*-algebras with a finite index \cite{KOT}. 

\begin{dfn}\label{Rohlin}
A conditional expectation $E$ of a unital  $C^*$-algebra $A$ with a finite index is said to have the {\it Rokhlin property} 
if there exists a  projection $e \in A_\infty$ satisfying 
$$
E^\infty(e) = ({\Index}E)^{-1} \cdot 1
$$
and a map $A \ni x \mapsto xe$ is injective. We call $e$ a Rokhlin projection.
\end{dfn}

\vskip 3mm

The following result says that 
the Rokhlin property of an action in the sense of Izumi implies that the canonical 
conditional expectation from a given C*-algebra to its fixed point algebra has the Rokhlin property in the sense of Definition~\ref{Rohlin}.

\vskip 1mm

\begin{prp}\label{finitegroup}(\cite{KOT})
Let $G$ be a finite group, $\alpha$  an action of $G$ on  a unital simple $C^*$-algebra $A$, and
 $E$  the canonical conditional expectation from $A$ onto the fixed point algebra $A^{\alpha}$.
Then $\alpha$ has the Rokhlin property if and only if $E$ has the Rokhlin property.
\end{prp}

\vskip 3mm

We need the following proposition to
prove the key lemma for the main theorem in the next section.

\begin{prp}\label{prp:tunnel}(\cite{KOT})
Let $P \subset A$ be an inclusion of unital $C^*$-algebras and $E$ a conditional expectation 
from $A$ onto $P$ with $\Index E < \infty $.
If there is a  projection $e \in A$ such that $E(e) = (\Index E)^{-1}$, 
then we have 
$$
ePe = Qe, \quad Q = P \cap \{e\}'
$$
In particular, if $e$ is a full projection, i.e., there are elements $x_i, y_i$ of $A$ such that 
$\sum_{i=1}^n x_i e y_i = 1$,  then 
$Q$ is a tunnel construction for $A \supset P$ such that $e$ is the Jones projection for $P \supset Q$.	
\end{prp} 

The following is a key lemma to prove the main theorem.

\begin{lem}\label{lem:embedding}
Let $P \subset A$ be an inclusion of unital \ca s and 
$E$  a conditional expectation from $A$ onto $P$ with a finite index.
If $E$ has the Rokhlin property with a Rokhlin projection $e \in A_\infty$, then 
there is a unital linear map 
$\beta \colon A^\infty 
\rightarrow P^\infty$ such that 
for any $x \in A^\infty$ there exists 
the unique element $y$ of $P^\infty$ such that $xe = ye = \beta(x)e$ and 
$\beta(A' \cap A^\infty) \subset P' \cap P^\infty$. 
In particular, $\beta_{|_A}$ is a unital injective *-homomorphism and 
$\beta(x) = x $ for all  $x \in P$.
\end{lem}

\begin{proof}
Let $e_P$  be the Jones projection for the inclusion $A \supset P$.
By the proof of Proposition \ref{prp:tunnel}, we have $ee_Pe = (\Index E)^{-1}e$. 
Therefore for any element $x$ in $A^\infty$
\begin{eqnarray*}
 xe&=& (\Index E)\hat{E}^\infty(e_P xe )\\
 &=& (\Index E)^2\hat{E}^\infty(e_P x e e_P e)\\
 &= &(\Index E)^2\hat{E}^\infty(E^\infty(xe)e_Pe) = (\Index E) E^\infty(xe)e,
\end{eqnarray*}
where $\hat{E}$ is the dual conditional expectation for $E$. Put $y = (\Index E) E^\infty(xe) \in P^\infty$.
 Then we have
$x e = y e$. 

Suppose that $ye = ze$ for $y, z \in P^\infty$. 
Then 
$$
z = z(\Index E) E^\infty(e) = (\Index E) E^\infty(ze) = (\Index E) E^\infty(ye) = y.
$$ 
Therefore we obtain the uniqueness of $y$.
Set $y = \beta(x)$.
Then $\beta$ is a unital linear map
from $A^\infty$ to $P^\infty$.
In particular, $\beta(x) = x$ for all $x \in P$.

\vskip 1mm

$\beta(A' \cap A^\infty) \subset P' \cap P^\infty$: $\forall x \in A' \cap A^\infty$ and $a \in P$ we have
\begin{align*}
\beta(x)ae &= \beta(x)ea\\
&= xea\\
&= xae\\
&= axe\\
&= a\beta(x)e
\end{align*}
From the uniqueness  we have
$\beta(x)a = a \beta(x)$, and 
$\beta(x) \in P' \cap P^\infty$.

\vskip 1mm

Homomorphism property of $\beta_{|_A}$: Let $x_1, x_2 \in A \subset A^\infty$ 
and $\beta(x_1), \beta(x_2) \in P^\infty$ 
such that $x_1e = \beta(x_1)e$, 
$x_2e = \beta(x_2)e$.

Then we have
\begin{align*}
x_1x_2e &=  x_1ex_2\\
&= \beta(x_1)ex_2\\
&= \beta(x_1)x_2e\\
&= \beta(x_1)\beta(x_2)e\\
\end{align*}
From the uniqueness we have
$\beta(x_1x_2) = \beta(x_1)\beta(x_2)$.

The *-preserving property of $\beta_{|A}$: 
Since for $x \in A$ $\beta_{|A}(x) = ({\Index}E)E^\infty(xe)$, 
\begin{align*}
\beta(x^*) &= ({\Index}E)E^\infty(x^*e)\\
&= ({\Index}E)E^\infty((xe)^*)\\
&= \{({\Index}E)E^\infty(xe)\}^*\\
&= \beta(x)^*
\end{align*}

Injectivity of $\beta_{|A}$:
Suppose that $\beta(x^*x) = 0$ for some $x \in A$, that is, 
$(\Index E) E^\infty(x^*xe) = 0$.
Since $e \in A' \cap A^\infty$, $E^\infty(x^*exe) = E^\infty((xe)^*(xe)) = 0$.
Since $E^\infty$ is faithful, we have $xe = 0$. Then from Definition~\ref{Rohlin} 
$x = 0$. This implies that $\beta$ is injective.

\end{proof}

\section{Rokhlin property and $\mathcal{D}$-absorbing}

Recall that a separable, unital C*-algebra $\mathcal{D}$ is called {\it strongly self-absorbing} 
if it is infinite-dimensional and the map $\id_{\mathcal{D}} \otimes 1_{\mathcal{D}}\colon
\mathcal{D} \rightarrow \mathcal{D} \otimes \mathcal{D}$ given by 
$d \mapsto d \otimes 1$ 
is approximately unitarily equivalent 
to an isomorphism $\varphi\colon \mathcal{D} \rightarrow \mathcal{D} \otimes \mathcal{D}$, 
that is, there is a suquence $(v_n)_{n\in\N}$ of unitaries in $\mathcal{D}$ 
satisfying 
$$
\|v_n^*(\id_{\mathcal{D}} \otimes 1_{\mathcal{D}}(d))v_n - \varphi(d)\|
\rightarrow 0 \ (n \rightarrow \infty) \ \forall d \in \mathcal{D}.
$$
A C*-algebra $A$ is called {\it $\mathcal{D}$-absorbing} if $A \otimes \mathcal{D} \cong A$. 
 
Note that if $A$ is an infinite dimensional unital simple AH algebra 
with no dimension growth, then $A$ is 
$\mathcal{Z}$-absorbing by \cite[Corollary 6.3]{Jiang-Su:absorbing}, where 
$\mathcal{Z}$ is called the Jiang-Su algebras, i.e., 
a direct limits of prime dimension drop algebras 
$I_{p,q} = \{f \in C([0, 1], M_{pq})\mid f(0) \in 1_p \otimes M_q, f(1) \in M_p \otimes 1_q\}$
for relative prime integers $p, q \geq 2$.
For the direct proof see \cite{DPT}.

In this section we consider the $\mathcal{D}$-absorbing property for 
an inclusion of unital C*-algebras. That is, 
let $P \subset A$ be an inclusion of unital C*-algebras and $E\colon A \rightarrow P$ 
be a conditional expectation of index finite with the Rokhlin property. 
If $A$ is $\mathcal{D}$-absorbing, then $P$ is $\mathcal{D}$-absorbing.

\vskip 3mm

We use the following characterization of the $\mathcal{D}$-absorbing.

\begin{thm}\label{thm:absorbing}(\cite[Theorem~7.2.2]{Ro})
Let $\mathcal{D}$ be a strongly self-absorbing and $A$ be any separable \ca. 
$A$ is $\mathcal{D}$-absorbing if and only if $\mathcal{D}$ admits a unital *-homomorphism to 
$A' \cap M(A)^\infty$.
\end{thm}

\vskip 3mm

The following is a special case of \cite[Lemma~2.4]{HW}. (See also \cite[Proposition~4.16]{HRW} 
and \cite[Remark~3.3.1]{W:Z-stable}.)

\vskip 1mm

\begin{lem}\label{lem:key2}
Let $A, B$ be unital separable C*-algebras. 
Suppose that for any finite sets $B_0 \subset B$, $A_0 \subset A$, and 
$\varepsilon > 0$ there is a completely positive contraction 
$\phi\colon B \rightarrow A^\infty$ such that for all $b, b' \in B_0$, 
$a \in A_0$ we have 

\begin{enumerate}
\item[(i)] $\|\phi(1) - 1\| < \varepsilon$,
\item[(ii)] $\|\phi(b)\phi(b') - \phi(bb')\| < \varepsilon$,
\item[(iii)] $\|[\phi(b), a]\| < \varepsilon$.
\end{enumerate}
Then there is a unital *-homomorphism from  $B$ to $A' \cap A^\infty$.
\end{lem}

\vskip 3mm

\begin{thm}\label{thm:main}
Let $P \subset A$ be an inclusion of unital 
C*-algebras and $E$ a conditional expectation from $A$ onto $P$ 
with a finite index.
Suppose that $\mathcal{D}$ is a separable unital self-absorbing 
\ca, $A$ is a separable $\mathcal{D}$-absorbing,
 and $E$ has the Rokhlin property. 
Then $P$ is $\mathcal{D}$-absorbing.
\end{thm}

\begin{proof}

Note that ${\mathcal D}$ is nuclear and simple by \cite[Lemma 3.10]{KP:embedding}.
Since $A$ is $\mathcal{D}$-absorbing, from Theorem~\ref{thm:absorbing} there is a unital *-homomorphism $\phi$ from $\mathcal{D}$ to $A' \cap A^\infty$.
Since $\mathcal{D}$ is nuclear, there is a completely positive lifting 
$\tilde{\phi}\colon {\mathcal D} \rightarrow \ell^\infty(A)$ by 
\cite{EH:lifting}, and 
write $\tilde{\phi} = (\phi_1, \phi_2, \dots, )$. 
Note that for all $d, d' \in {\mathcal D}$ and $a \in A$
\begin{enumerate}
\item[(i)]
$$
\lim\sup_i\|\phi_i(d)\phi_i(d') - \phi_i(dd')\| = 0,
$$
\item[(ii)]
$$
\lim\sup_i\|[\phi_i(d), a]\| = 0.
$$
\end{enumerate}

Since $E$ has the Rokhlin property,
there is a unital *-homomorphism $\beta$ from $A$ to $P^\infty$ 
such that $\beta(x) = x$ for all $x \in P$ by Lemma~\ref{lem:embedding}.
Set $\psi_i = \beta \circ \phi_i$. Then 
$\psi_i$ is a unital completely positive map from $\mathcal{D}$ to 
$P^\infty$ such that for all $b, b' \in {\mathcal D}$ and $a \in P$ 

\begin{enumerate}
\item[(iii)] 
\begin{align*}
\|\psi_i(d)\psi_i(d') - \psi_i(dd')\| 
&= \|\beta \circ \phi_i(d)\beta \circ \phi_i(d') - \beta \circ \phi_i(dd')\|\\
&= \|\beta (\phi_i(d)\phi_i(d') - \phi_i(dd'))\|\\
&\leq \|\phi(d)\phi(d') - \phi_i(dd')\|
\end{align*}
\item[(iv)]
\begin{align*}
\|[\psi_i(d), a]\| &= \|[\beta \circ \phi_i(d), \beta(a)]\| \ (\beta(a) = a)\\
&= \|\beta(\phi_i(d))\beta(a) - \beta(a)\beta(\phi_i(d))\|\\
&= \|\beta(\phi_i(d)a - a\phi_i(d)\|\\
&\leq \|[\phi_i(d), a]\|.
\end{align*}
\end{enumerate}

Let ${\mathcal D_0} \subset {\mathcal D}$ and  $P_0 \subset P$ 
be any finite sets.
For any $\varepsilon > 0$  taking sufficient large $k$
then 
from (i) $\sim$ (iv) we have for all $d, d' \in {\mathcal D_0}$
and $a \in P_0$

\begin{enumerate}
\item[(v)] $\|\psi_k(d)\psi(d') - \psi_k(bb')\| < \varepsilon$,
\item[(vi)] $\|[\psi_k(d), a]\| < \varepsilon$.
\end{enumerate}

Hence we have a unital *-homomorphism from ${\mathcal D}$ to 
$P' \cap P^\infty$ by Lemma~\ref{lem:key2}. 

Therefore $P$ is $\mathcal{D}$-absorbing from Therem~\ref{thm:absorbing}.
\end{proof}

\vskip 3mm

\begin{cor}\label{Cor:Crossedproducts}(\cite{HW})
Let $A$ be a separable, unital, simple $\mathcal{D}$-absorbing C*-algebra and 
$\alpha$ be  an action of a finite group $G$ on $A$.
Suppose that $\alpha$ has the Rohklin property. Then 
the crossed product algebra $A \rtimes_\alpha G$ is $\mathcal{D}$-absorbing.
\end{cor}

\begin{proof}
Since $\alpha$ is outer, the canonical conditional expectation 
$E\colon A \rightarrow A^\alpha$ is of index finite type, where 
$E(a) = \frac{1}{|G|}\sum_{g\in G}\alpha_g(a)$, 
where $|G|$ is the order of $G$.
Then $E$ has the Rokhlin property by \cite[Proposition 2.9]{KOT}
 and $A^\alpha$ is $\mathcal{D}$-stable by Theorem~\ref{thm:main}.
 
 Consider the following basic construction:
 $$
 A^\alpha \subset A \subset B_1 (= A \rtimes_\alpha G).
 $$
 Note that $B_1$ is isomorphic to $pM_n(A^\alpha)p$ for some $n \in \N$
 and a projection $p \in M_n(A^\alpha)$ by \cite[Lemma~3.3.4]{Watatani:index}. 
 Since $M_n(A^\alpha)$ is $\mathcal{D}$-absorbing and the $\mathcal{D}$-absorbing
 is stable under the hereditary by \cite[Corollary 3.1]{TW}, we know that 
 $B_1$ is also $\mathcal{D}$-absorbing, that is, 
 $A \rtimes_\alpha G$ is $\mathcal{D}$-absorbing.
\end{proof}

\vskip 3mm

Recall that a separable C*-algebra $A$ is said to be {\it approximately divisible} 
if there is a finite dimensional C*-algebras $B$ with no abelian summands, 
which admits a unital embedding into $A' \cap M(A)^\infty$, where $M(A)$ means the 
multiplier algebra of $A$. 
In \cite{TW2} Toms and Winter proved that separable, approximately divisible
C*-algebras are $\mathcal{Z}$-absorbing. 
We show that separable, approximately divisibility  is stable under inclusion of 
index finite with the condition that  the correspondent  conditonal expactation has the Rokhlin property. 

\vskip 3mm

\begin{prp}\label{prp:app divisible}
Let $P \subset A$ be an inclusion of unital C*-algebra and 
$E$ a conditional expectation from $A$ onto $P$ with a finite index.
Suppose that $A$ is a separable, approximately divisible and $E$ has the Rokhlin property. 
Then $P$ is approximately divisible.
\end{prp}

\vskip 3mm

\begin{proof}
Since $A$ is approximately divisible, there is a finite dimensional 
C*-algebra $B$ with no abelian summands which admits 
a unital embedding $\phi$ from $B$ into $A' \cap A^\infty$. 
Let $\beta\colon A \rightarrow P^\infty$ be a unital injective *-homomorphism
in Lemma~\ref{lem:embedding}. 
Since $B$ is nuclear and finite direct sums of simple C*-algebras, 
through the same argument in Theorem~\ref{thm:main}
we have a unital embedding from $B$ into $P' \cap P^\infty$, 
hence $P$ is approximately divisible.
\end{proof}

\vskip 3mm

\section{Rokhlin property and strongly self-absorbing}

Recall that separable unital C*-algebra $\mathcal{D}$ is said to 
have {\it approximately inner half flip} if 
the two natural inclusions of $\mathcal{D}$ into $\mathcal{D} \otimes \mathcal{D}$
as the first and second factor, repectively, are approximately unitarily equivalent, i.e.,
there is a sequence $(v_n)_{n\in \N}$ of unitaries in $\mathcal{D} \otimes \mathcal{D}$  such that $\|v_n(d \otimes 1_D)v_n^* - 1_D \otimes d\| \rightarrow 0 (n \rightarrow \infty)$ 
for $d \in \mathcal{D}$. 
In \cite{ER:inner flip} Effros and Rosenberg proved that if  $A$ is AF C*-algebra, $A$ has approximate  
half inner flip if and only if $A$ is a UHF algebra.
Note that a separable unital C*-algebra $A$ has approximately inner half-flip 
implies that $A$ is simple and nuclear by \cite[Propositions~2.7 and 2.8]{ER:inner flip}.

Under the condition that separable unital C*-algebra $\mathcal{D}$ 
has approximately inner half flip, Toms and Winter gave the characterization when 
$\mathcal{D}$ is strongly self-absorbing \cite{TW}.
Using this characterization we show that if a conditional expectation 
$E\colon \mathcal{D} \rightarrow P$ 
for an inclusion of separable unital C*-algebras $P \subset \mathcal{D}$ with index finite,
has the Rokhlin property and  $\mathcal{D}$  is semiprojective and strongly self-absorbing, 
then $P$ is strongly self-absorbing.

\vskip 3mm

\begin{prp}\label{Prp:half flip}
Let $P \subset A$ be an inclusion of separable unital C*-algebras 
with index finite and $A$ have approximately inner half flip. 
Suppose that $E$ has the Rokhlin property and $A$ is semiprojective.
Then $P$ has approximately inner half flip.
\end{prp}

\vskip 3mm

\begin{proof}
We have only to show that for any finite set $\mathcal{F} \subset P$ 
and $\varepsilon > 0$ there exists an unitary $w \in P \otimes P$
such that

$$
\|w(a \otimes 1_P)w^* - 1_P \otimes a\| < \varepsilon \quad \forall a \in \mathcal{F}.
$$

Let $\mathcal{F} \subset P$ with $\|a\| \leq 1$ for any $a \in \mathcal{F}$ 
and $\varepsilon > 0$.
Since $A$ has approximately inner half flip, there exists 
a unitary $u \in A \otimes A$ such that 
$$
\|u(a \otimes 1_P)u^* - 1_P \otimes a\| < \frac{1}{5}\varepsilon \quad \forall a \in \mathcal{F}.
$$
Note that $1_A = 1_P$.

Let $\beta\colon A \rightarrow P^\infty$ be a *-homomorphism in 
Lemma~\ref{lem:embedding} such that $\beta(a) = a$ for $a \in P$.
Since $A$ is semiprojective, there exists $n \in \N$ and 
*-homomorphism $\hat{\beta}\colon A \rightarrow \ell^\infty(\N, P)/I_n$ 
such that $\beta = \pi \circ \hat{\beta}$, where 
$I_n = \{(a_k)_{k=1}^\infty \in \ell^\infty(P)\colon a_k = 0\ \hbox{for} \ k > n\}$
and $\pi\colon \ell^\infty(\N, P)/I_n \rightarrow P^\infty$ is a canonical quotient map. 
Note that $P^\infty = \ell^\infty(\N, P)/\overline{\{\cup_{n=1}^\infty I_n\}}$.
Since $\beta(a) = a$ for $a \in P$, if we write 
$\hat{\beta} = (\phi_1, \phi_2, \dots, \phi_n, \dots)$, then 
$$
\lim\sup\|a - \phi_k(a)\| = 0.
$$
Hence there is a *-homomorphism $\phi_k$ for $k > n$ such that 
$$
\|\phi_k(a) - a\| < \frac{1}{5}\varepsilon
$$
for $a \in \mathcal{F} \cup \{1_P\}$.
We have then for any $a \in \mathcal{F}$ 
\begin{align*}
&\|(\phi_k\otimes\phi_k)(u)(a \otimes 1_P)(\phi_k\otimes\phi_k)(u^*) - 1_P \otimes a\| \\
&= \|(\phi_k\otimes\phi_k)(u)((a - \phi_k(a) + \phi_k(a))\otimes 1_P)(\phi_k\otimes\phi_k)(u^*) - 1_P \otimes a\|\\
&\leq 
\|(\phi_k\otimes\phi_k)(u)((a - \phi_k(a))\otimes 1_P)(\phi_k\otimes\phi_k)(u^*)\| \\
&+ \|(\phi_k\otimes\phi_k)(u)(\phi_k(a))\otimes (1_P - \phi_k(1_P) + \phi_k(1_P))(\phi_k\otimes\phi_k)(u^*) - 1_P \otimes a\|\\
&\leq 
\|a - \phi_k(a)\| 
+ \|(\phi_k\otimes\phi_k)(u)(\phi_k(a))\otimes (1_P - \phi_k(1_P))(\phi_k\otimes\phi_k)(u^*)\|\\
&+ 
\|(\phi_k\otimes\phi_k)(u)(\phi_k(a))\otimes \phi_k(1_P)(\phi_k\otimes\phi_k)(u^*) - 1_P \otimes a\|\\
&\leq \|a - \phi_k(a)\| + \|1_P - \phi_k(1_P)\| + 
\|(\phi_k\otimes\phi_k)(u)(\phi_k(a))\otimes \phi_k(1_P)(\phi_k\otimes\phi_k)(u^*) - 
\phi_k(1_P) \otimes \phi_k(a)\| \\
&+ \|\phi_k(1_P - \phi_k(1)) \otimes a + \phi_k(a) \otimes (a - \phi_k(a))\|\\
&\leq 2\|a - \phi_k(a)\| + 2\|1_P - \phi_k(1_P)\|
 + \|(\phi_k \otimes \phi_k)(u(a \otimes 1_P)u^* - 1_P \otimes a)\|\\
&< 5 \times \frac{1}{5}\varepsilon = \varepsilon.
\end{align*}
Since $(\phi_k \otimes \phi_k)(u) \in P \otimes P$ is unitary, 
we conclude that $P$ has approximately inner half flip map.
\end{proof}

\vskip 5mm 

\begin{rmk}
Under the same condition for an inclusion of separable unital C*-algebras
$P \subset A$ in Theorem~\ref{Prp:half flip} 
since $P$ has approximately inner half flip map
we know that $P$ is nuclear and simple.
But when there exists a conditional expectation $E\colon A \rightarrow P$, 
which is called a projection of norm one, if $A$ is nuclear,
then we can easily prove that $P$ is simple. Indeed, since $A$ is nuclar, for any finite set
$\mathcal{F} \subset P \subset A$ and $\varepsilon > 0$ there exists a finite rank completely 
positive contructive $\phi\colon A \rightarrow A$ such that 
$$
\|\phi(a) - a\| < \varepsilon \quad \forall a \in \mathcal{F}
$$
from \cite[Theorem~3.1]{CE:nuclear}. 
Since $E$ is completey positive \cite[Theorem~3.4]{Ta:book1}, 
$E \circ \phi_{|P} \colon P \rightarrow P$ 
is finite rank completely positive contructive, and for any $a \in \mathcal{F}$ 
we have  
\begin{align*}
\|(E \circ \phi_{|P})(a) - a\| 
&= \|E(\phi(a) - a)\|\\
&\leq \|\phi(a) - a\| < \varepsilon.
\end{align*}
Hence $P$ is nuclear.

To deduce the simplicity of $P$ we need the Rokhlin condition for $E\colon A \rightarrow P$ 
as follows. 

Let $P \subset A$ be an inclusion of separable unital C*-algebras 
with index finite.
Suppose that $A$ is simple and $E$ has the Rokhlin property.
 
Let $e$ be the Rokhlin projection for $E$. Then the dual conditional expectation 
$\hat{E}\colon A_1 \rightarrow A$ has the approximately representable
by \cite[Proposition~3.4 (1)]{KOT}, 
that is, we have  
$exe = \hat{E}(x)e$ for any $x \in A_1$, where $P \subset A \subset A_1$ 
is the basic construction. 
Then $A' \cap A_1 \subset A$ by \cite[Proposition~3.5]{KOT}. 
Since $A$ is simple, we have

\begin{align*}
A_1' \cap A_1 &\subset A' \cap A_1\\
&\subset A' \cap A \cong \C.
\end{align*}
 
Then $A_1$ is simple by \cite[Theorem~3.3]{Izumi:inclusion}. 
Since $P$ is stably isomorphic to $A_1$, we conclude that $P$ is simple.
\end{rmk}

\vskip 5mm

If $\mathcal{D}$ is a strongly self-absorbing inductive limit of recursive
subhomogeneous algebras in the sense of Phillips \cite{Phillips:recursive}, then 
$\mathcal{D}$ is either projectionless (i.e. the Jiang-Su algebra $\mathcal{Z}$) 
or a UHF algebra of infinite type  by \cite[Corollary~5.10]{TW}.
On the contrary, if $\mathcal{D}$ is a separable purely infinite strongly self-absorbing 
C*-algebra which satisfies the Universal Coefficients Theorem 
(We write $\mathcal{D}$ is in the UCT class $N$.). Then 
$\mathcal{D}$ is either $\mathcal{O}_2$, $\mathcal{O}_\infty$ or 
a tensor product of $\mathcal{O}_\infty$ with a UHF algebra of inifinite type by 
\cite[Corollary~5.2]{TW}.

At first we show that the universal UHF algebra $\mathcal{U}_\infty$ is stable 
under crossed products by finite groups actions with the Rokhlin property.
However, this is not true in the case of an action with the tracial Rokhlin property 
which is weeker than the Rokhlin property.
(See section 6.)

\vskip 3mm

For $n \in \N$ $M_{n^\infty}$ means the UHF algebra of type $n^\infty$.

\vskip 3mm

The following lemma may be well known result, but we put the proof for a self-contained.

\vskip 3mm

\begin{lem}\label{lem:UHF}
Let $A$ be a UHF algebra and $\alpha$ be an action from a finite group $G$
on $A$. Suppose that $\alpha$ has the Rokhlin property. Then
$A \rtimes_\alpha G$ is a UHF algebra.
\end{lem}

\vskip 3mm

\begin{proof}
Since $\alpha$ has the Rokhlin property, for any finite set $F \subset A \rtimes_\alpha G$
and $\varepsilon > 0$, there are $n$, projection $f \in A$, and  unital homomorphism
$\varphi\colon M_n \otimes fAf \rightarrow A\rtimes_\alpha G$ such that 
$\dist(a, \varphi(M_n \otimes fAf)) < \varepsilon$ for all $a \in F$ by \cite[Theorem~3.2]{OP:Rohlin}.
Since $A$ is a UHF algbera, $M_n \otimes fAf$ is also a UHF algebra. Hence we may 
assume that there is a full matrix subalgebra $M_k$ of $A \rtimes_\alpha G$ 
such that $\dist(a, M_k) < \varepsilon$.
This implies that $A \rtimes_\alpha G$ is a UHF algebra by \cite[Therem~1.13]{Glimm:certain}.
\end{proof}

\vskip 3mm

\begin{thm}\label{thm:strongly absorbing}
Let $\mathcal{D}$ be $\mathcal{U}_\infty$ and 
let $\alpha$ be an action of a finite group $G$ on $\mathcal{D}$.
Suppose that $\alpha$ has the Rokhlin property. Then 
the crossed prodct $\mathcal{U}_\infty\rtimes_\alpha G$ is isomorphic to $\mathcal{U}_\infty$.
\end{thm}

\vskip 3mm

\begin{proof}
Since $\mathcal{U}_\infty$ is in the UCT class $\mathcal{N}$, 
has tracial topological rank zero in the sense of Lin, 
and $\alpha$ has the Rokhlin property with $\alpha_* = 0$ on $K_*(\mathcal{U}_\infty)$ 
for $i = 1, 2$, $(\mathcal{U}_\infty, \alpha)$ is conjugate to 
$(\mathcal{U}_\infty \otimes M_{|G|^\infty}, \id_{\mathcal{U}_\infty} \otimes \mu^G)$ 
by \cite[Theorem~3.5]{Izumi:Rohlin2}, 
that is, there is an isomorphism $\theta\colon \mathcal{U}_\infty \rightarrow 
\mathcal{U}_\infty \otimes M_{|G|^\infty}$ such that 
$$
\theta\circ \alpha_g\circ \theta^{-1} = (\id_{\mathcal{U}_\infty} \otimes \mu^G)_g \ 
\forall g \in G,
$$
where $\lambda$ is the left regular representation of $G$ and 
$\mu_g^G = \otimes_{n=1}^\infty Ad(\lambda(g))$ \ $(g \in G )$.
Note that $B(\ell^2(G))$ is identified with $M_{|G|}$ and $\mu^G$ has the 
Rokhlin property.

Then 
$\mathcal{U}_\infty \rtimes_\alpha G$ is isomorphic to 
$(\mathcal{U}_\infty \otimes M_{|G|^\infty})
\rtimes_{\id_{\mathcal{U}_\infty} \otimes \mu^G} G$.
Since 
\begin{align*}
(\mathcal{U}_\infty \otimes M_{|G|^\infty})
\rtimes_{\id_{\mathcal{U}_\infty} \otimes \mu^G} G
&= \mathcal{U}_\infty \otimes (M_{|G|^\infty}\rtimes_{\mu^G}G)\\
&\cong \mathcal{U}_\infty,
\end{align*}
since $M_{|G|^\infty}\rtimes_{\mu^G}G$ is a UHF algebra by Lemma~\ref{lem:UHF}.

Therefore $\mathcal{U}_\infty \rtimes_\alpha G$ is strongly self-absorbing.
\end{proof}

\vskip 5mm

The following example implies that 
the Rokhlin property is essential in Theorem~\ref{thm:strongly absorbing}.

\vskip 3mm

\begin{exa}
Let $\mathcal{U}_\infty$ be the universal UHF algebra and 
$A = M_{2^\infty}$. Then $A \otimes \mathcal{U}_\infty \cong \mathcal{U}_\infty$.

Let $\alpha$ be an symmetry in \cite[Proposition~5.1.2]{B:symmetry}. 
Then $A \rtimes_\alpha \Z/2\Z$ is not a AF algebra.
We note that $\alpha$ has the tracial Rokhlin property 
by \cite[Proposition~3.4]{Phillips:cyclic}, but does not 
have the Rokhlin property, since the crossed product algebra 
$A \rtimes_\alpha \Z/2\Z$ is not AF algebra by \cite[Theorem~2.2]{Phillips:tracial}. 

Then $\alpha \otimes id$ is a symmetry  with the tracial Rokhlin property 
on $A \otimes \mathcal{U}_\infty (\cong A)$, and the crossed product algebra 
\begin{align*}
(A \otimes \mathcal{U}_\infty) \rtimes_{\alpha \otimes id}\Z/2\Z 
&\cong (A \rtimes_\alpha\Z/2\Z) \otimes \mathcal{U}_\infty\\
&\cong B \otimes \mathcal{U}_\infty,
\end{align*}
where $B$ is the Bunce-Dedens algebras of type $2^\infty$ by \cite[Proposition~5.4.1]{B:symmetry}.
Note that $K_1(B \otimes \mathcal{U}_\infty) \not= 0$, that is, 
$B \otimes \mathcal{U}_\infty$ is not a AF algebra.
Since a strongly self-absorbing inductive limit of type I with real rank zero 
C*-algebra is a UHF algebra of infinite type by \cite[Corollary~5.9]{TW},  
$B \otimes \mathcal{U}_\infty$ is not a strongly self-absorbing algebra.
Hence there is a symmetry $\beta$ with the tracial Rokhlin property on $\mathcal{U}_\infty$ such that 
$\mathcal{U}_\infty \rtimes_\beta \Z/2\Z$ is not strongly self-absorbing.
\end{exa}

\vskip 5mm

Next we show  that strongly self-absorbing property is stable 
under the assumption that a conditional expectation $E$ for an inclusion of unital C*-algebras 
$P \subset A$ with index finite has the Rokhlin property  and $A$ is semiprojective.
Note that Cuntz algebras $O_2$ and $O_\infty$ are semiprojective (\cite{B:Shape}\cite{B:semiprojective}).

\vskip 5mm

\begin{thm}\label{thm:strongly self-absorbing}
Let $P \subset A$ be an inclusion of unital separable C*-algebras with index finite.
Suppose that a conditional expectation $E\colon A \rightarrow P$ has the 
Rokhlin property and $A$ is semiprojective and strongly self-absorbing.
Then $P$ is strongly self-absorbing.
\end{thm}

\vskip 3mm

\begin{proof}
Since $A$ is strongly self-absorbing, $A$ has  approximately inner half flip 
by \cite[Proposition~1.5]{TW}. Hence we know that $P$ has also approximately 
inner half flip by Proposition~\ref{Prp:half flip}.

To prove the strongly self-absorbing for $P$ we have only to show that 
there are a unital *-homomorphism $\gamma\colon P \otimes P \rightarrow P$ and an 
approximately central sequence of unital endmorphims of $P$ by \cite[Proposition~1.10(ii)]{TW}.

Since $A$ is strongly self-absorbing, there are a unital *-homomorphism 
$\gamma\colon A \otimes A \rightarrow A$ and for any finite set 
$\mathcal{F}$ and $\mathcal{G}$ in $A$ and $\varepsilon > 0$ 
there exists a unital endmorphism $\psi\colon A \rightarrow A$ such that 
$$
\|[\psi(a), b]\| < \varepsilon
$$
for $a \in \mathcal{F}$ and $b \in \mathcal{G}$.

Consider a unital *-hmomorphism $\beta$ from $A$ into $P^\infty$ 
in Lemma~\ref{lem:embedding}.
Since $A$ is semiprojective, there are a $n \in \N$ and a unital 
*-homomorphism $\overline{\beta}\colon A \rightarrow \ell^\infty(P)/I_n$ 
such that $\pi \circ \overline{\beta} = \beta$, where
$I_n = \{(a_k)_{k=1}^\infty \in \ell^\infty(\N, P)\colon a_k = 0\ \hbox{for} \ k > n\}$
and $\pi\colon \ell^\infty(P)/I_n \rightarrow P^\infty$ is a canonical quotient map.
Let write $\overline{\beta} = (\phi_1,\phi_2,\dots)$, then for $k > n$ $\phi_k$ 
is a unital *-homomorphism.
Then $\phi_k \circ \gamma_{| P\otimes P}\colon P \otimes P \rightarrow P$ 
is a unital injective *-homomorphism. Note that $P$ is simple.

Let $\mathcal{F}$ and $\mathcal{G}$ in $P$ be finite sets 
with $\|a\| , \|b\| \leq 1$ for $a \in \mathcal{F}$ and $b \in \mathcal{G}$, 
and $\varepsilon > 0$.
Take $\psi\colon A \rightarrow A$ such that 
$$
\|[\psi(a), b]\| < \frac{1}{4}\varepsilon
$$
for $a \in \mathcal{F}$ and $b \in \mathcal{G}$. 
Hence 
\begin{align*}
\|(\phi_k (\psi(a)b - b\psi(a) )\| &< \frac{1}{4}\varepsilon\\
\|(\phi_k\circ\psi)(a)\phi_k(b) - \phi_k(b)(\phi_k\circ\psi)(a)\| & < \frac{1}{4}\varepsilon.
\end{align*}

Since $\lim\sup_k\|\phi_k(a) - a \| = \|\beta(a) - a\| = 0$ for $a \in P$, there is a $k \in \N$ 
such that 
$$
\|\phi_k(b) - b\| < \frac{1}{4}\varepsilon
$$
for $b \in {\mathcal{G}}$.
Then we have 
\begin{align*}
\|[\phi_k\circ \psi(a), b]\| &= \|[\phi_k\circ \psi(a),\phi_k(b) + b - \phi_k(b)] \| \\
&\leq \|[\phi_k\circ \psi(a),\phi_k(b)]\| + \|[\phi_k\circ \psi(a),b -\phi_k(b)]\|
&< \frac{3}{4}\varepsilon < \varepsilon.
\end{align*}
Hence there exist an approximately central sequence of unital *-homomorphisms of $P$.

Hence $P$ is strongly self-absorbing by \cite[Proposition~1.10(ii)]{TW}. 
\end{proof}

\vskip 5mm

\begin{cor}\label{cor:self-absorbing}
Let $P \subset A$ be an inclusion of unital separable C*-algebras with index finite.
Suppose that a conditional expectation $E\colon A \rightarrow P$ has the 
Rokhlin property. Suppose that $A$ is 
$O_2$ or $O_\infty$.  Then $P \cong A$.
\end{cor}

\vskip 3mm

\begin{proof}
Since $O_2$ and $O_\infty$ are separable, semiprojective, and strongly self-absorbing,
$P$ is also strongly self-absorbing by Theorem~\ref{thm:strongly self-absorbing}.
As in the proof of Theorem~\ref{thm:strongly self-absorbing} there exist unital embeddings
$\iota_P\colon P \rightarrow A$ and $\iota_A\colon A \rightarrow P$.
Hence we conclude that $P \cong A$ by  \cite[Proposition~5.12]{TW}.
\end{proof}

\vskip 5mm

\begin{cor}\label{cor:Cuntz}(\cite[Theorem~4.2]{Izumi:Rohlin1})
Let $\alpha$ be an action of a finite group $G$ on $\mathcal{O}_2$. 
Suppose that $\alpha$ has the Rohklin property.
Then we have 
\begin{enumerate}
\item
$\mathcal{O}_2^G \cong \mathcal{O}_2$.
\item 
The crossed product algebra $O_2 \rtimes_\alpha G \cong O_2$.
\end{enumerate}
\end{cor}

\vskip 3mm

\begin{proof}
Since $\alpha$ is outer, the canonical conditional expectation 
$E\colon \mathcal{O}_2 \rightarrow \mathcal{O}_2^G$ is of index finite type, where 
$E(a) = \frac{1}{|G|}\sum_{g\in G}\alpha_g(a)$, 
where $|G|$ is the order of $G$.

$(1)$: Since $E$ has the Rokhlin property by \cite[Proposition 2.9]{KOT}, 
$\mathcal{O}_2^G \cong \mathcal{O}_2$ by Theorem~\ref{cor:self-absorbing}.

$(2)$: Consider the following basic construction:
 $$
 \mathcal{O}_2^\alpha \subset \mathcal{O}_2 \subset B_1 (= \mathcal{O}_2 \rtimes_\alpha G).
 $$
 Note that $B_1$ is isomorphic to $pM_n(\mathcal{O}_2^\alpha)p$ for some $n \in \N$
 and a projection $p \in M_n(\mathcal{O}_2^\alpha)$ by \cite[Lemma~3.3.4]{Watatani:index}. 
 Since $M_n(\mathcal{O}_2^\alpha)$ is isomorphic to $\mathcal{O}_2$ and $\mathcal{O}_2$
 is stable under the hereditary, we know that 
 $B_1$ is isomorphic to $\mathcal{O}_2$, that is, 
$\mathcal{O}_2 \rtimes_\alpha G$ is isomorphic to $\mathcal{O}_2$.
\end{proof}

\vskip 5mm

\begin{rmk}
From \cite[Theorem~3.6]{Izumi:Rohlin2} there is no non-trivial finite group action with the Rokhlin property
on $\mathcal{O}_\infty$
\end{rmk}

\vskip 5mm


\section{Local $\mathcal{C}$-property}

Through this section $\mathcal{D}$ means the separable unital strongly 
self-absorbing C*-algebra.

We recall the definition of local $\mathcal{C}$-property in 
\cite{OP:Rohlin}.

\begin{dfn}\label{D:FSat}
Let ${\mathcal{C}}$ be a class of separable unital C*-algebras.
Then ${\mathcal{C}}$ is {\emph{finitely saturated}}
if the following closure conditions hold:
\begin{enumerate}
\item\label{D:FSat:1}
If $A \in {\mathcal{C}}$ and $B \cong A,$ then $B \in {\mathcal{C}}.$
\item\label{D:FSat:2}
If $A_1, A_2, \ldots, A_n \in {\mathcal{C}}$ then
$\bigoplus_{k=1}^n A_k \in {\mathcal{C}}.$
\item\label{D:FSat:3}
If $A \in {\mathcal{C}}$ and $n \in \N,$
then $M_n (A) \in {\mathcal{C}}.$
\item\label{D:FSat:4}
If $A \in {\mathcal{C}}$ and $p \in A$ is a nonzero projection,
then $p A p \in {\mathcal{C}}.$
\end{enumerate}
Moreover,
the {\emph{finite saturation}} of a class ${\mathcal{C}}$ is the
smallest finitely saturated class which contains ${\mathcal{C}}.$
\end{dfn}

\vskip 3mm

\begin{rmk}
Let $\mathcal{C}$ be the set of all separable,  unital, 
$\mathcal{D}$-absorbing C*-algebras.
Then ${\mathcal{C}}$ is finitely saturated.
\end{rmk}

\begin{proof} It is obvious that $\mathcal{C}$ satisfies 
three conditions $(1) \sim (3)$.
The condition $(4)$ comes from \cite[Corollary~3.1]{TW}.
\end{proof}

\vskip 3mm

\begin{dfn}\label{D:LC}
Let ${\mathcal{C}}$ be a class of separable unital C*-algebras.
A {\emph{unital local ${\mathcal{C}}$-algebra}}
is a separable unital C*-algebra $A$
such that for every finite set $S \subset A$ and every $\ep > 0,$
there is a C*-algebra $B$ in the finite saturation of ${\mathcal{C}}$
and a unital *-homomorphism $\ph \colon B \to A$
(not necessarily injective)
such that $\dist (a, \, \ph (B)) < \ep$ for all $a \in S.$
\end{dfn}

\vskip 3mm

\begin{thm}\label{Th:LocalZ-absorbing}
Let $\mathcal{C}$ be the class of 
all separable, unital, $\mathcal{D}$-absorbing C*-algebras.
Let $A$ be a local $\mathcal{C}$-algebra. Then $A \in \mathcal{C}$.
\end{thm}

\vskip 3mm

\begin{proof}
Since $A$ is a local $\mathcal{C}$-algebra,
for any finite sets $D_0 \subset D$ and $F \subset A$, and 
$\varepsilon > 0$   there exists $B_0 \in \mathcal{C}$ 
and a unital homomorphisms $\psi \colon B_0 \rightarrow A$ such that 
$\dist(f, \psi(B_0)) < \frac{\varepsilon}{4}$ for any $f \in F$. 
That is, there exists $b_f \in B_0$ such that 
$\|f - \psi(b_f)\| < \varepsilon$ for each $f \in F$. We may assume that 
the norm of each element in $D_0$ is less than or equalt to one.

Since $B_0$ is $\mathcal{D}$-absorbing, there is a complete positive contractive 
$\phi\colon D \rightarrow B_0$ such that 

\begin{enumerate}
\item
$\|\phi(1) - 1\| < \varepsilon$,\\
\item 
$\|\phi(d)\phi(d') - \phi(dd')\| < \varepsilon$,\\
\item
$\|[\phi(d), b]\| < \frac{\varepsilon}{2}
\ \forall d, d' \in D_0, \forall b \in \{b_f\}_{f \in F}$ 
\end{enumerate}
(See the first part of the proof in Theorem~\ref{thm:main}.)

\vskip 1mm

Then $\psi \circ \phi \colon D \rightarrow A$ is a complete positive 
contractive satisfying 
\begin{align*}
\|\psi\circ\phi(1) - 1\| 
&\leq  \|\phi(1) - 1\| < \varepsilon \\
\|\psi\circ\phi(d)\psi\circ\phi(d') - \psi\circ\phi(dd')\| 
&= \|\psi(\phi(d)\phi(d') - \phi(dd'))\|\\
&\leq \|\phi(d)\phi(d') - \phi(dd')\| < \varepsilon\\
\|[\psi\circ\phi(d), f]\| &= \|[(\psi\circ\phi(d), \psi(b_f)]\| 
+ 2\|f - \psi(b_f)\|\|\psi\circ\phi(d)\|\\
&< \frac{\varepsilon}{2} + 2 \times \frac{\varepsilon}{4}\times 1 
= \varepsilon
\end{align*}
for any $d, d' \in D_0$ and $f \in F$.

\vskip 2mm

Hence for any finite sets $D_0 \subset D$ and $F \subset A$, and 
$\varepsilon > 0$ 
there is a completely positive contractive map 
$\varphi\colon D \rightarrow A \hookrightarrow A^\infty$ such that 

\begin{enumerate}
\item
$\|\varphi(1) - 1\| < \varepsilon$,\\
\item 
$\|\varphi(d)\varphi(d') - \varphi(dd')\| < \varepsilon$,\\
\item
$\|[\varphi(d), f]\| < \varepsilon, \ \ 
\ \forall d, d' \in D_0,\  \forall f \in F$
\end{enumerate}

Therefore, there is a unital homomorphism 
$\rho\colon D \rightarrow A' \cap A^\infty$ by Lemm~\ref{lem:key2}, and 
$A$ is $\mathcal{D}$-absorbing by Theorem~\ref{thm:absorbing}. 
\end{proof}

\vskip 3mm

We will give another proof of Corollary~\ref{Cor:Crossedproducts}:

\begin{cor}\label{Cor:AnotherCrossedProducts}(Corollary~\ref{Cor:Crossedproducts})
Let $A$ be a separable unital, $\mathcal{D}$-absorbing C*-algebra and 
$\alpha$ be  an action of a finite group $G$ on $A$.
Suppose that $\alpha$ has the Rokhlin property. Then 
the crossed product algebra $A \rtimes_\alpha G$ is $\mathcal{D}$-absorbing.
\end{cor}

\vskip 3mm

\begin{proof}
For any $\varepsilon > 0$ and a finite set $F \subset A \rtimes_\alpha G$ 
there are $n \in \N$, 
a projection $f \in A$, and a unital homomorphism 
$\varphi\colon M_n \otimes fAf \rightarrow A \rtimes_\alpha G$ such that 
$\dist(a, \varphi(M_n \otimes fAf)) < \varepsilon $ for all $a \in F$ by 
\cite[Theorem 3.2]{OP:Rohlin}. 

Since $fAf$ is $\mathcal{D}$-absorbing by \cite[Corollary 3.1]{TW}, 
$M_n \otimes fAf$ is also $\mathcal{D}$-absorbing. We also know that 
$\varphi(M_n \otimes fAf)$ is $\mathcal{D}$-absorbing by 
\cite[Corollary 3.3]{TW}. Hence $A \rtimes_\alpha G$ is local $\mathcal{C}$-algebra 
in the same class $\mathcal{C}$ in Theorem~\ref{Th:LocalZ-absorbing}.
Therefore, $A \rtimes_\alpha G$ is $\mathcal{D}$-absorbing.
\end{proof}

\vskip 3mm

\begin{rmk}
From the previous sequence works by  Kodaka and the authos \cite{KOT}, and the first author and 
Phillips \cite{OP:Rohlin} we can conclude the followings.

Let $P \subset A$ be an inclusion of separable unital C*-algebras and 
$E\colon A \rightarrow P$ be of index finite type with the Rokhlin property. 
Suppose that $A$ belongs to a class $\mathcal{C}$ characterized by some algebraic structural 
property. Then $P$ belongs to the class $\mathcal{C}$. 
The classes we consider include:

\begin{enumerate}
\item
C*-algebras with various kind of direct limit decompositiions involving semiprojective building blocks.
(For examples, AF algebras, AI algebras, AT algebras etc.)
\item
Simple unital AH algebras with slow dimension growth and real rank zero.
\item
C*-algebras with real rank zero.
\item
C*-algebras with stable rank one.
\item
Simple C*-algebras for which the order on projections is determined by traces.
\item 
C*-algebras with isometrically extremal richness.
\item
C*-algebras with approximately divisivity.
\item
C*-algebras with ${\mathcal D}$-absorbing for a separable unital strongly self-absorbing 
C*-algebra ${\mathcal D}$.
\end{enumerate}
\end{rmk}

\vskip 3mm

\section{Intermediate fixed point algebras}

In this section we present an inclusion of unital C*-algebras 
$P \subset A$ which does not come from an action of finite group on 
$A$.

Let $A$ be a unital C*-algebra and $\alpha$ an action of a finite group on $A$.
Suppose that $\alpha$ has the Rokhlin property. Then $\alpha$ is outer, that is, 
for any $g \in G\backslash\{1\}$ $\alpha_g$ is outer. 
Hence $\alpha$ is saturated, that is, the set 
$\{\hat{x}\hat{y}^*\colon x, y \in  A \}$ spans 
a dense subalgebra in $A \rtimes_\alpha G$, where $\hat{x} = \sum_{g\in G}\alpha_g(x)\lambda_g$ 
for $x \in A$.
Then the canonical conditional expectation $E \rightarrow A^G$ is of 
index finite type with index $|G|$ by \cite[Theorem~4.1]{Jeong:saturated} and has 
the Rokhlin property by \cite[Proposition~2.9]{KOT}. 
In particular, by \cite[Theorem~4.1]{Jeong:saturated} 
we have a quasi-basis for $\{(u_i, u_i^*)\}$ for $E$ 
such that 
\begin{enumerate}
\item for any $x \in A$
$$
x = \sum_{i=1}^nE(xu_i)u_i^*.
$$
\item
$$
{\rm Index}E = \sum_iu_iu_i^* = |G|\\
$$
\item
$$
\sum_iu_i\alpha_g(u_i^*) = 0, \ g \not= 1
$$
\end{enumerate}

Moreover, for any 
subgroup $H $ of $G$ the correspondent conditional expectation 
$F = E_{|H}\colon A^H \rightarrow A^G$ is of index finite and has the Rokhlin property. 

Indeed, consider a set $\{E_H(u_i), E_H(u_i^*)\}$ for $A^H \times A^H$, 
where $E_H\colon A \rightarrow A^H$ is the canonical conditional expectation by 
$$
E_H(x) = \frac{1}{|H|}\sum_{h\in H}\alpha_h(x).
$$
Then $\{E_H(u_i), E_H(u_i^*)\}$ becomes a quasi-basis for $F$ as follows:
for any $x \in A^H$
\begin{align*}
\sum_iF(xE_H(u_i))E_H(u_i^*) &= \sum_iE(xE_H(u_i))E_H(u_i^*)\\
&= \sum_iE_H(E(xE_H(u_i))u_i^*)\\
&= \sum_iE_H(\frac{1}{|G|}\sum_{g\in G}\alpha_g(xE_H(u_i))u_i^*)\\
&= \frac{1}{|G|}\frac{1}{|H|}E_H(\sum_i(\sum_{g \in G}\alpha_g(x)\alpha_g(\sum_{h\in H}\alpha_h(u_i))u_i^*)\\
&= \frac{1}{|G|}\frac{1}{|H|}E_H(\sum_{h\in H}\sum_{g\in G}\alpha_g(x)\sum_i\alpha_{gh}(u_i)u_i^*)\\
&= \frac{1}{|G|}\frac{1}{|H|}E_H(\sum_{h \in H}\alpha_{h^{-1}}(x)\sum_iu_iu_i^*)\ ((3))\\
&= \frac{1}{|H|}E_H(\sum_{h\in H}\alpha_h(x))\\
&= E_H(E_H(x)) = x
\end{align*}
We also have 
\begin{align*}
{\rm Index}F &= \sum_iE_H(u_i)E_H(u_i^*)\\
&= \sum_i\frac{1}{|H|}\sum_{h\in H}\alpha_h(u_i)(\frac{1}{|H|}\sum_{k\in H}\alpha_k(u_i^*))\\
&= \frac{1}{|H|^2}\sum_i\sum_{h\in H}\sum_{k\in H}\alpha_h(u_i)\alpha_k(u_i^*)\\
&= \frac{1}{|H|^2}\sum_{h\in H}\sum_{k\in H}\alpha_h(\sum_iu_i\alpha_{h^{-1}k}(u_i^*))\ ((3))\\
&= \frac{1}{|H|^2}\sum_{h\in H}\alpha_h(\sum_iu_iu_i^*)\\
&= \frac{1}{|H|^2}\times |H| \times |G| = \frac{|G|}{|H|}
\end{align*}

To prove the Rokhlin property we set $e_H = \sum_{h\in H}e_h$. 
Since $A' \subset (A^H)'$, $e_H \in (A^H)' \cap (A^H)^\infty$. 
We have then
\begin{align*}
F^\infty(e_H) &= \sum_{h\in H}F^\infty(e_h)\\
&= \frac{1}{|G|}\sum_{g\in G}\alpha_g(\sum_{h \in H}e_h)\\
&= \frac{1}{|H|}\sum_{h\in H}\sum_{g\in G}\alpha_{g}(e_h)\\
&= \frac{|G|}{|H|}1
\end{align*}
Therefore $e_H$ is the Rokhlin projection for $F$.

\vskip 5mm

\begin{prp}\label{prp:intermidiate H}
In the same condition in the above we could conclude that 
$E_H\colon A \rightarrow A^H$ has the Rokhlin property for any 
subgruoup $H$ of $G$.
\end{prp}

\vskip 3mm

\begin{proof}
Let $G = \cup_{i=1}^lHg_i$ be a decomposition of right cosets in $G$
and $f = \sum_{i=1}^le_{g_i}  \in A' \cap A^\infty$ be a projection.
Then $f$ is  the Rokhlin projection for $E_H$. Indeed 

\begin{align*}
E_H^\infty(f) &= \sum_{i=1}^lE_H^\infty(e_{g_i})\\
&= \frac{1}{|H|}\sum_{i=1}^l\sum_{h\in H}\alpha_h^\infty(e_{g_i})\\
&= \frac{1}{|H|}\sum_{i=1}^l\sum_{h\in H}e_{hg_i}\\
&= \frac{1}{|H|}\sum_{g\in G}e_g\\
&= \frac{1}{|H|}
\end{align*}

Since index of $E_H = \frac{1}{|H|}$, we conclude that
$E_H$ has the Rokhlin property.
\end{proof}

\vskip 5mm

From the above observation we have the following characterization.

\begin{prp}\label{prp:intermediate}
Let $A$ be a separable unital C*-algebra, $\alpha$ an action of 
a finite group $G$ on $A$ and $E\colon A \rightarrow A^G$ a canonical 
conditional expectation.
Suppose that $\alpha$ has the Rokhlin property. Then we have 
\begin{enumerate}
\item 
For any subgroup $H$ of $G$ 
the restricted $E$ to $A^H$, which is a conditional expectation from $A^H$ onto $A^G$, has the 
Rokhlin property.
\item
Let $\mathcal{D}$ be a separable unital strongly self-absorbing C*-algebra 
and $A$ be $\mathcal{D}$-absorbing.
Then for any subgroup $H$ of $G$ $A^H$ is $\mathcal{D}$-absorbing.
\item
If $A = \mathcal{O}_2$, then for any subgroup $H$ of $G$ $A^H \cong \mathcal{O}_2$.
\end{enumerate}
\end{prp}

\vskip 3mm

\begin{proof}
We use the same notations in the above argument.

$(1)$ follows from the above argument.

$(2)$: 
From Proposition~\ref{prp:intermidiate H} 
a conditional expectation $E_H\colon A \rightarrow A^H$ 
is of index finite with the Rokhlin proeprty. 

Since $A$ is $\mathcal{D}$-absorbing, so is $A^H$ by 
Theorem~\ref{thm:main}. 

$(3)$: Since $E_H\colon A \rightarrow A^H$ has the Rokhlin property by 
Proposition~\ref{prp:intermidiate H}, $A^H \cong \mathcal{O}_2$ 
by Corollary~\ref{cor:self-absorbing}. 
\end{proof}

\vskip 5mm

\begin{rmk}
Let $A$ be a unital C*-algebra and $\alpha$ be an action from a finite group $G$ on $A$.
Let $H$ be a subgroup of $G$. Then the condition that an inclusion $A^G \subset A^H$ is isomorphic to 
$B^K \subset B$ for some C*-algebra $B$ and an action from a finite group $K$ on $B$ 
implies that $H$ is a normal subgroup of $G$ (c.f. \cite{T:normal intermediate}).
Hence from Proposition~\ref{prp:intermediate} we have examples of conditional expectations 
for inclusions of unital C*-algebras with the Rokhlin property which do not come from 
finite group actions.
\end{rmk}


\end{document}